\documentstyle[12pt]{article}
\advance\hoffset by -0.5 truecm

\newtheorem{Theorem}{Theorem}[section]
\newtheorem{Lemma}[Theorem]{Lemma}
\newtheorem{Corollary}[Theorem]{Corollary}

\def\df{\buildrel \rm def \over =} 
\newtheorem{Proposition}[Theorem]{Proposition}

\def \dim{{\hbox {\rm dim}}}
\def \Vol{{\hbox {\rm Vol}}}
\def\df{\buildrel \rm def \over =} 
\def\PROOF{{\em Proof}\,: }              
\def\QED{~\hfill~ $\diamond$ \vspace{7mm}}

\def\ep{\varepsilon}

\begin{document}
\title{Counting geodesics on a Riemannian manifold and topological entropy of geodesic flows}
\author{Keith Burns\thanks{Partially supported by NSF grant DMS-9206923}
\and
Gabriel P. Paternain}

\date{June 1995}
\maketitle

\centerline{\em To Ricardo Ma\~n\'e, in Memoriam}

\begin{abstract}Let $M$ be a compact $C^{\infty}$ Riemannian manifold.
Given $p$ and $q$ in $M$ and $T>0$, define $n_{T}(p,q)$ as the number of 
geodesic segments joining $p$ and $q$ with length $\leq T$.
Ma\~n\'e showed in \cite{M} that 
\[\lim_{T\rightarrow \infty}\frac{1}{T}\log \int_{M\times M}n_{T}(p,q)\,dpdq
 = h_{top},\]
where $h_{top}$ denotes the topological entropy of the geodesic flow of $M$.

In the present paper we exhibit an open set of metrics on the two-sphere for which

\[\limsup_{T\rightarrow\infty}\frac{1}{T}\log n_{T}(p,q)< h_{top},\] 
for a positive measure set of $(p,q)\in M\times M$.
This answers in the negative questions raised by Ma\~n\'e in \cite{M}.

\end{abstract}

\section{\em Introduction}

Let $M$ be a closed connected 
 $C^{\infty}$ manifold and $g$ a $C^{r}$ ($r\geq 2$) Riemannian metric on $M$. 
Let $\phi_{t}:SM\rightarrow SM$ denote the geodesic flow of $g$, acting on the 
unit sphere bundle $SM$. Denote the topological entropy of the 
flow $\phi_{t}$ with respect to a compact subset $K$ of $SM$ by $h_{top}(K)$. 
Then $h_{top} \df h_{top}(SM)$ is the topological entropy of $\phi_t$. 

Given $p$ and $q$ in $M$ and $T>0$, define $n_{T}(p,q)$ as the number of 
geodesic segments joining $p$ and $q$ with length $\leq T$. Already in 1962
Berger and Bott \cite{BB} observed that there are significant relationships
between integrals of this function and the dynamics of the geodesic flow. As
was pointed out in \cite{BB}, it is not hard to 
see that, for each $T>0$, the counting function $n_{T}(p,q)$ is finite and locally constant on an 
open full measure subset of $M\times M$, and integrable on $M \times M$. More generally, if $N$ is a compact submanifold of $M$, define $n_T(N,q)$ to be the number of geodesic segments with length $\leq T$ that join a point in $N$ to $q$ and are initially orthogonal to $N$. The function $n_T(N,q)$ enjoys properties similar to those of $n_T(p,q)$.

For $M$ and $N$ that are $C^\infty$, it was shown in \cite{PP} that Yomdin's theorem \cite{Y} can be used to prove that
\begin{equation}
\limsup_{T\rightarrow \infty}\frac{1}{T}\log \int_{M}n_{T}(N,q)\,dq
\leq h_{top}(S^\perp N),  
\label{PPineq}
\end{equation}
where $S^\perp N$ is the set of unit vectors with footpoint in $N$ that are orthogonal to $N$. It was also shown in \cite{PP} that when $N$ is the diagonal in the product manifold $M \times M$, the above inequality reduces to  
\begin{equation}
\limsup_{T\rightarrow \infty}\frac{1}{T}\log \int_{M\times M}n_{T}(p,q)\,dpdq
\leq h_{top} . 
\label{2}
\end{equation}
Another proof of (\ref{2}) is given in \cite{M}.
 On the other hand, for any $C^{r}$ Riemannian 
metric ($r\geq 2$), Ma\~{n}\'{e} shows in \cite{M} that
 \begin{equation}
\liminf_{T\rightarrow\infty}\frac{1}{T} \log \int_{M\times M}n_{T}(p,q)\,dpdq
\geq h_{top} . 
\label{3}
\end{equation}

Ma\~{n}\'{e} thereby  obtains the first purely Riemannian characterization of 
$h_{top}$ for an arbitrary $C^{\infty}$ Riemannian metric:  combining 
(\ref{2}) and (\ref{3}) gives

\begin{equation}
\lim_{T\rightarrow \infty}\frac{1}{T}\log \int_{M\times M}n_{T}(p,q)\,dpdq
 = h_{top}.  
\label{4}
\end{equation}

Ma\~n\'e's result extends  earlier work of Manning and Freire-Ma\~n\'e.
Suppose $\widetilde p$ is a lift of a point $p \in M$ to the universal cover 
$\widetilde M$ of $M$ and $B_T(\widetilde p)$ is the ball of radius $T$ about 
$\widetilde p$ in $\widetilde M$ (with the metric lifted from $M$). It follows 
easily from the results in \cite{BB} that
\begin{equation}
 \int_M n_{T}(p,q)\,dq \geq  \Vol\,B_T(\widetilde p),  
\label{4.5}
\end{equation}
with equality if $M$ has no conjugate points. Manning showed that, 
for any $\widetilde p \in \widetilde M$, \ $T^{-1} \Vol\,B_T(\widetilde p)$
converges to a limit $\lambda$ that is independent of $\widetilde p$. From 
(\ref{4}) and (\ref{4.5}), it is easy to obtain the inequality $h_{top} \geq 
\lambda$ for any Riemannian manifold, which was  first proved by Manning in \cite{Ma}.
One also sees that $h_{top} = \lambda$ if $M$ has no conjugate points, which was first proved by Freire and Ma\~n\'e in \cite{FM}.

Besides the natural appeal of a formula like (\ref{4}), 
there are other reasons to be interested in relations between the topological 
entropy of the geodesic flow and the growth rate of the 
average number of geodesic segments between two points in the manifold. The 
function $n_{T}(p,q)$ also counts the number of critical points of the energy 
functional on the path space $\Omega^{T}(p,q)$ given by all the curves joining 
$p$ and $q$ with length $\leq T$. Therefore using Morse theory, $n_{T}(p,q)$ 
can be bounded from below by the sum of the Betti numbers of $\Omega^{T}(p,q)$ 
(provided of course that $p$ and $q$ are not conjugate).  By averaging over 
$M$ and using results of Gromov \cite{Gr2}, one can obtain in this fashion 
remarkable relations between the topology of $M$ and the vanishing of 
$h_{top}$; we refer to \cite{Gr1,P,P1,PP} for a detailed description of these 
relations.

The present paper continues the study of relationships between $h_{top}$ and
the exponential growth rate of $n_T(p,q)$. The strongest relationship is  
for metrics with no conjugate points: Ma\~n\'e showed in \cite{M} that in this
case
\[ \lim_{T \to \infty}\frac1T \log n_T(p,q) = h_{top} \quad\hbox{\rm  for all $p,q \in M$.}\]
 For any $C^\infty$ Riemannian manifold, taking $N$ to be the single point $p$ in (\ref{PPineq}) gives the following inequality, that was first proved by G.P. Paternain  in \cite{P}:
 \begin{equation}
 \limsup_{T\rightarrow\infty} \frac{1}{T} \log\int_{M}n_{T}(p,q)\,dq
 \leq h_{top}(S_pM)  \leq h_{top}\quad\hbox{\rm for all $p \in M$.}  \label{5}
 \end{equation} 
 As Ma\~{n}\'{e} indicates in \cite{M}, a simple argument using the Borel-Cantelli Lemma shows 
that for every $p \in M$ one has
\begin{equation}
 \limsup_{T \to \infty}\frac{1}{T}\log n_{T}(p,q) \leq
\limsup_{T\rightarrow\infty} \frac{1}{T} \log\int_{M}n_{T}(p,q')\,dq'
\quad\hbox{\rm for a.e. $q \in M$.} \label{5.1}
\end{equation}
It is immediate from  (\ref{5}) and (\ref{5.1}) that  for all $p\in M$ one has
 \begin{equation}
\limsup_{T\rightarrow\infty}\frac{1}{T}\log n_{T}(p,q)\leq h_{top} 
\quad \hbox{for a.e. $q\in M$.}    \label{5.2}
\end{equation}
This inequality motivated Ma\~{n}\'{e} to pose the following questions in 
\cite{M}. \bigskip
 
\noindent{\bf Question I.} {\em Is it true that
 \begin{equation}
\lim_{T\rightarrow\infty}\frac{1}{T}\log n_{T}(p,q)= h_{top} 
\quad\hbox{\rm  for a.e. $(p,q)\in M\times M$?} 
\label{6}
\end{equation}}

\noindent{\bf Question II.} {\em Is it true that equation (\ref{6}) holds for generic Riemannian 
metrics when~$\dim\, M = 2$?}

\medskip

The main purpose of the present paper is to give a negative answer to Question II.
Of course this also answers Question I negatively.

Because of inequalities (\ref{5}) and (\ref{5.1}) it is natural to consider modified versions of the first question. \bigskip

\noindent{\bf Question I$'$.} {\em Is it true that
 \[\limsup_{T\rightarrow\infty}\frac{1}{T}\log\int_{M} n_{T}(p,q)\,dq 
  =   h_{top} \quad\hbox{ for all $p \in M$?}\]}

\noindent{\bf Question I$''$.} {\em Is it true that 
 \[\limsup_{T\rightarrow\infty}\frac{1}{T}\log\int_{M} n_{T}(p,q)\,dq 
  =   h_{top} \quad\hbox{for almost every $p \in M$?}\]}

We shall also give negative answers to Questions I$'$ and I$''$.
Because of (\ref{5.1}), a negative answer to Question I$''$ implies a negative answer to Question I. Our main example will be a surface for which Question I$''$ has a negative answer. The important features of this surface will be stable under small perturbations of the metric; we thereby obtain an open set of metrics for which the answer to Question I$''$ is negative. This implies a negative answer to Question II.

Let us describe the contents of the paper in more detail. 
 In Section 2 we obtain a slight improvement of inequality (\ref{5}) that is 
needed for our main example.
  We also prove an interesting related result, namely, that for 
the case of surfaces, the set of points $p\in M$ for which inequality (\ref{5}) is an equality has non-empty interior. 
This result will imply an alternative proof of Ma\~{n}\'{e}'s inequality (\ref{3}) for the two dimensional case.

In Section 3, we observe that the Weinstein examples, described in \cite{BeBe}, give rise to many manifolds 
 (such as ${\bf CP}^{k}$) 
 for which there exists a point $p$ such that $\displaystyle\int_{M}n_{T}(p,q)\,dq$ grows {\bf linearly} in $T$, even though $h_{top}>0$.
  It follows from (\ref{5.1}) that $\limsup_{T \to \infty} T^{-1} \log n_T(p,q) = 0$
for a.e. $q \in M$. 
 All these examples have dimension 
$\geq 3$. We also construct a metric on the two-sphere with $h_{top}>0$ which has a point $p$ such that all 
the geodesics leaving from $p$ are simple, closed and have the same period, so again $\displaystyle\int_{M}n_{T}(p,q)\,dq$ grows {\bf linearly} in $T$.

Section 4 contains our main example that gives negative answers to 
Ma\~{n}\'{e}'s Questions I and II and also to Question I$''$.
 We start with a one parameter family of 
surfaces of revolution $C_{d}$, $0< d<\infty$, which have  the 
properties shown in Figure 1.

\begin{figure} 

\vspace{8cm}

\caption{Main example.}

\end{figure}

Each $C_d$ is diffeomorphic to the two sphere and contains a region $R$ bounded by a geodesic circle of latitude $\alpha$ which is shorter than any other geodesic circle of latitude in $R$.
Inside $R$ is a circle of latitude $\gamma_0$ that is a hyperbolic closed geodesic. The region $R$ contains geodesics that are both forward and backwards asymptotic to  $\gamma_{0}$. 
This means that $\gamma_{0}$ has a homoclinic connection.
Attached to $R$ there is a flat cylinder of length $d$ which we close smoothly with a cap $D$. Let $P$ be the center of the cap and $Q$ the center of the region $R$.

We perturb the metric in a small part of $R$ (shaded in Figure 1) so as to break the homoclinic 
connection of $\gamma_{0}$ and obtain a horseshoe with entropy $h > 0$.
 Moreover we arrange the perturbation so that the meridians of $C_{d}$ are preserved, i.e., all geodesics leaving from $P$ in the new metric are simple, closed, have the same period and they coincide outside $R$ with the meridians of the old metric.
 A careful application of KAM theory (Lemma \ref{key}) shows that  there is  
$d_0 > 0$ such that for any $d\geq d_0$ the geodesic flow for the perturbed metric has an invariant torus whose projection to $C_{d}$ becomes singular on two curves, one near $P$, the other near $Q$, as shown in Figure 1.

For each $d$, the torus and its image under the flip $v\rightarrow -v$ separate the unit sphere bundle of $C_{d}$ into two invariant 
sets, $W_{1}$ and $W_{2}$.
This forces the geodesics which pass sufficiently close to $P$ or $Q$ to cross 
the piece of the surface around $\gamma_{0}$ fairly perpendicular to 
$\gamma_{0}$; intuitively this separating surface keeps the geodesics that leave from points sufficiently close to $P$ or $Q$ away from the horseshoe. This property allows us to 
estimate their Liapunov exponents and make them small compared with $h$, by 
choosing $d$ large enough. 
  If $W_{2}$ is the set that contains the horseshoe, then the topological 
entropy on $W_{2}$ is at least $h$ and the topological entropy on $W_{1}$ can 
be made smaller than $h$, since the Lyapunov exponents in $W_1$ can be made small, as we explained above. The results in Section 2 will imply 
that if $p$ is sufficiently close to $P$ or $Q$ , then 

\[
\limsup_{T\rightarrow\infty}\frac{1}{T}\log n_{T}(p,q)
 \leq 
\limsup_{T\rightarrow\infty}\frac{1}{T}\log\int_{M}n_{T}(p,q')\,dq'\leq h_{top}(W_{1}), \]
for almost every $q\in M$. 
Thus there exists a positive measure set $U\subset M\times M$ such that for 
$(p,q)\in U$ 
\begin{equation}
\limsup_{T\rightarrow\infty}\frac{1}{T}\log n_{T}(p,q)\leq 
\limsup_{T\rightarrow\infty}\frac{1}{T}\log\int_{M}n_{T}(p,q)\,dq'<h \leq 
h_{top}.
\label{Upq}
\end{equation} 
 This gives negative answers to Question I and I$''$.
 A careful look at the details in 
Section 4 will show that for small perturbations of this example  
(\ref{Upq}) still holds for a positive measure set of $p$ and $q$. Thus Question II also has a negative answer.

{\sl Acknowledgements:} 
The second author is grateful to Detlef Gromoll for suggesting several years ago the study of the relationship between the exponential growth rate of $n_{T}(p,q)$ and $h_{top}$.
He is also grateful to the University of Maryland and Northwestern University for their hospitality while this work began. 
The first author thanks the Universidad de la Rep\'ublica, Uruguay, for hospitality while this work was completed.

\section{\em Some properties of the function $n_{T}(p,q)$}

We begin this section by proving a result (cf. Proposition \ref{p1} below) that we shall use in our main example. In what follows we will always assume that the Riemannian metric is $C^{\infty}$.

In \cite{M}, Ma\~{n}\'{e} shows the following application of the Borel-Cantelli Lemma:

\begin{Lemma}Let $(X,{\cal A},\mu)$ be a probability space and $f_{n}:X\rightarrow (0,+\infty)$ a sequence of integrable functions. Then:
\[\limsup_{n\rightarrow\infty}\frac{1}{n}\log f_{n}(x)\leq\limsup_{n\rightarrow\infty}\frac{1}{n}\log\int_{X}f_{n}\,d\mu,\]
for $\mu$-a.e. $x\in X$.  \label{Borel}
\end{Lemma}

As in the Introduction, denote by $SM$ the unit sphere bundle of $M$ and let $\pi:SM\rightarrow M$ be the canonical projection.
For $p\in M$, set $S_{p}= \pi^{-1}(p)$.
If $K\subset SM$ is a closed set, let $h_{top}(K)$ denote the topological entropy of the geodesic flow $\phi_{t}$ with respect to the set $K$; with this notation, $h_{top}(SM)=h_{top}$.
Combining Lemma (\ref{Borel}) above with Corollary 2.2 from \cite{PP} gives the following result; we include the proof for the convenience of the reader.

\begin{Proposition}Let $K$ be any closed set in $SM$ and suppose there exists $p\in M$, such that $S_{p}\subset K$, then
\[\limsup_{T\rightarrow\infty}\frac{1}{T}\log n_{T}(p,q)\leq 
\limsup_{T\rightarrow\infty}\frac{1}{T}\log\int_{M}n_{T}(p,q')\,dq'\leq h_{top}(K),\]
for almost every $q\in M$.     \label{p1}
\end{Proposition}

\PROOF It is known (cf. \cite{BB,P,PP}) that for any $p\in M$,
\begin{equation}
\int_{M}n_{T}(p,q)\,dq\leq\int_{0}^{T}\Vol(\phi_{t}S_{p})dt,   \label{7'}
\end{equation}
where ``Vol" stands for the $n-1$ dimensional Riemannian volume ($n= \dim\, M$).
Now we use Yomdin's Theorem as stated in \cite{Gr1} to obtain:
\[\limsup_{T\rightarrow\infty}\frac{1}{T}\log \Vol(\phi_{T}S_{p})\leq h_{top}(S_{p}).\]
Since $S_{p}\subset K$, $h_{top}(S_{p})\leq h_{top}(K)$ and by combining the last inequality with (\ref{7'}) we obtain
\[\limsup_{T\rightarrow\infty}\frac{1}{T}\log\int_{M}n_{T}(p,q)\,dq\leq h_{top}(K),\]
and by using Lemma \ref{Borel} we conclude:
\[\limsup_{T\rightarrow\infty}\frac{1}{T}\log n_{T}(p,q)\leq 
\limsup_{T\rightarrow\infty}\frac{1}{T}\log\int_{M}n_{T}(p,q')\,dq'\leq h_{top}(K),\]
for almost every $q\in M$.
\QED

For brevity, let $\sigma_{p}$ be:
 \[\sigma_{p}\df\limsup_{T\rightarrow\infty}\frac{1}{T}\log\int_{M}n_{T}(p,q)\,dq.\]
 Then (\ref{5}) can be writen as $\sigma_{p}\leq h_{top}$ for all $p\in M$.
Our main example in Section 4 possesses the following property: there exists an open set $V\subset M$ such that for any $p\in V$, $\sigma_{p}<h_{top}$.
To complete the description of the possible behavior of the
correspondence $p\rightarrow \sigma_{p}$ when $\dim\,M=2$, we now show that in this case 
 the set of points for which $\sigma_{p}=h_{top}$ always has
non-empty interior.

\begin{Theorem}Suppose $\dim\,M = 2$. Consider the set $\Omega$ of
points $p\in M$ such that
\[\lim_{T\rightarrow\infty}\frac{1}{T}log\int_{M}n_{T}(p,q)\,dq=h_{top}.\]
Then $\Omega$ has non-empty interior.  \label{teo}
\end{Theorem}

\PROOF For each vector $v\in S_{p}$ consider the Jacobi equation along
the geodesic $\gamma_v$ defined by $v$:
\begin{equation}
y''(t)+K(t)y(t)=0,   \label{jacobi}
\end{equation}
where $K(t)$ is the Gaussian curvature of $M$ at $\gamma_v(t)$. 
Let $y_{v}(t)$ be the solution of
(\ref{jacobi}) that satisfies: $y_{v}(0)=0$ and $y'_{v}(0)=1$.
On account of the results in \cite{BB} we have: 
\begin{equation}
\int_{M}n_{T}(p,q)\,dq=\int_{0}^{T}ds\int_{S_{p}} | y_{v}(s) | \,dv.
\label{one}
\end{equation}
On the other hand by definition, if $l$ denotes length, we can write
\begin{equation}
l(\phi_{T}S_{p})=\int_{S_{p}}\sqrt{y_{v}^{2}(T)+(y_{v}')^{2}(T)}\,dv.
\label{two}
\end{equation}
Equations (\ref{one}) and (\ref{two}) clearly imply
\begin{equation}
\limsup_{T\rightarrow\infty}\frac{1}{T}\log \int_{M}n_{T}(p,q)\,dq\leq\limsup_{T\to\infty}\frac{1}{T}\log l(\phi_{T}S_{p}).
  \label{three}
\end{equation}

Since $y'_{v}(0)=1$, it follows from (\ref{jacobi}) that
\[y_{v}'(T)=1-\int_{0}^{T}K(s)y_{v}(s)\,ds;\]
 if $L$ is the maximum of the absolute value of the curvature on $M$, we obtain
\[ |  y_{v}'(T) | \leq 1+L\int_{0}^{T} |  y_{v}(s) |  \,ds.\]
Therefore using (\ref{two}) we can write
\begin{eqnarray*}
l(\phi_{T}S_{p}) &\leq& \int_{S_{p}} \left( |  y_{v}(T) | 
+1+L\int_{0}^{T} |  y_{v}(s) |  \,ds \right) \,dv  \\
&=& \int_{S_{p}} |  y_{v}(T) | \,dv+l(S_{p})+L\int_{0}^{T}\int_{S_{p}} 
| y_{v}(s) |  \,dvds. \end{eqnarray*}
{}From the last inequality and (\ref{one}) we obtain
\begin{equation}
\liminf_{T\rightarrow\infty}\frac{1}{T}\log l(\phi_{T}S_{p})\leq \liminf_{T\rightarrow\infty}\frac{1}{T}\log \int_{M}n_{T}(p,q)\,dq.  \label{four}
\end{equation}

Recall now from \cite[page 220]{N}, that there exists an arc $\alpha\subset SM$ such that
\begin{equation}
\lim_{T\rightarrow\infty}\frac{1}{T}\log l(\phi_{T}\alpha)=h_{top}.  \label{six}
\end{equation}
Moreover, if $\mu$ is an ergodic measure of maximal entropy, $\alpha$ could be {\sl any} arc transversal to the Pesin stable manifolds of $\mu$.
Let $v_{0}$ be a Pesin point of $\mu$ and let $W^{s}(v_{0})$ denote the weak stable manifold through $v_{0}$.
 By a well known property of the geodesic flow ($W^{s}(v_{0})$ is $\phi_{t}$-invariant), there exists $r>0$ such that $T_{\phi_{r}v_{0}}W^{s}(v_{0})\cap T_{\phi_{r}v_{0}}S_{\pi(\phi_{r}v_{0})}=\{0\}$, therefore the curve $S_{\pi(\phi_{r}v_{0})}$ is transversal to $W^{s}(v_{0})$ at the point $\phi_{r}v_{0}$.
 Hence (\ref{three}), (\ref{four}) and (\ref{six}) imply:
\begin{equation}
\lim_{T\rightarrow\infty}\frac{1}{T}\log \int_{M}n_{T}(p,q)\,dq=h_{top},                                                              \label{seven}
\end{equation}
where $p=\pi(\phi_{r}v_{0})$.
 On the other hand, $\pi:W^{s}(v_{0})\rightarrow M$ is a local diffeomorphism at $\phi_{r}v_{0}$, since $T_{\phi_{r}v_{0}}W^{s}(v_{0})\cap T_{\phi_{r}v_{0}}S_{\pi(\phi_{r}v_{0})}=\{0\}$.
Therefore every unit circle $S_{p'}$, with foot point $p'$ in a neighborhood of $p=\pi(\phi_{r}v_{0})$, is transversal to $W^{s}(v_{0})$ and thus (\ref{seven}) has to hold for an open set around $p$.
\QED

We finish this section by showing how Theorem \ref{teo} implies Ma\~{n}\'{e}'s inequality (\ref{3}) for the two dimensional case. 

\begin{Corollary}Suppose $\dim\,M = 2$. Then:
\[\liminf_{T\rightarrow\infty}\frac{1}{T} \log \int_{M\times M}n_{T}(p,q)\,dpdq\geq h_{top}.\]
\end{Corollary}

\PROOF For $p\in M$, set $I_{p} = \int_{M}n_{T}(p,q)\,dq$.
Let $\Omega$ be the set from Theorem \ref{teo} and let $m$ denote its
measure (by Theorem \ref{teo}, $m>0$).
Then
\[ \int_{M\times M}n_{T}(p,q)\,dpdq=\int_{M}I_{p}(T)\,dp\geq \int_{\Omega}I_{p}(T)\,dp,\]
and by Jensen's inequality
\[\log \int_{M\times M}n_{T}(p,q)\,dpdq\geq\frac{1}{m}\int_{\Omega}\log I_{p}(T)\,dp.\]
Hence
\[\liminf_{T\rightarrow\infty}\frac{1}{T} \log \int_{M\times M}n_{T}(p,q)\,dpdq\geq\frac{1}{m}\int_{\Omega}\sigma_{p}\,dp=h_{top}.\]
\QED

\section{\em Simple examples}

We consider in this section the following modified version of Question I:\medskip

\noindent{\bf Question I$'$.} {\em Is it true that
\[\limsup_{T\rightarrow\infty}\frac{1}{T}\log \int_{M}n_{T}(p,q)\,dq=h_{top}
\quad\hbox{ for all $p$?}\]}

Suppose $M^{n}=D\cup_{S^{n-1}} D_{N}$, where $D$ is an $n$-dimensional disk with center $p$ and $D_{N}$ is a disk bundle over a manifold $N$, such that its associated sphere bundle $\partial D_{N}$ is diffeomorphic to $S^{n-1}$.
 A typical example of such a manifold is ${\bf CP}^{k}$; if we remove a disk from it, we obtain a disk bundle over ${\bf CP}^{k-1}$.
Let $g_{N}$ be any Riemannian metric on $N$, then it is shown in \cite{BeBe} that one can construct a metric $g$ on $M$ so that:
\begin{itemize}
\item all geodesics leaving from $p$ (the center of $D$) return to $p$ exactly at the same time.
\item $(N,g_{N})\rightarrow (M,g)$ is a totally geodesic isometric embedding.
\end{itemize}

Choose  a metric $g_{N}$ for which the topological entropy of the geodesic flow of $N$ is positive. Since $N$ is totally geodesic, this implies that $h_{top}$ of the geodesic flow of $M$ is also positive.
On the other hand since every geodesic leaving from $p$ returns to $p$ at exactly the same time, it follows that $\Vol(\phi_{T}S_{p})$ is uniformly bounded for all $T$ and by equation (\ref{7'}), $\displaystyle\int_{M}n_{T}(p,q)\,dq$ grows linearly with $T$.
  This gives a negative answer to Question I$'$.
In fact, it also shows that for some $p$ the growth of $\displaystyle\int_{M}n_{T}(p,q)\,dq$ could be only linear, even though $h_{top}>0$.

Observe that this construction requires $\dim\,N \geq 2$, and hence $\dim\,M \geq 3$, so it is natural to ask if similar examples can be constructed in the case of surfaces.
We shall show below that there exist metrics on $S^{2}$ with $h_{top}>0$ and the additional property that there exists a point $p$ such that all the geodesics leaving from $p$ are simple, closed and with the same period. 
As a consequence, $\displaystyle\int_{M}n_{T}(p,q)\,dq$ grows {\sl linearly} with $T$.

Let $C$ be a surface of revolution diffeomorphic to $S^{2}$.
On $C$ the geodesic flow is completely integrable with an integral of motion --- the Clairaut integral --- given by $r\sin\phi$, where $r$ is the radial distance of a point on the surface to the axis of revolution and $\phi$ is the angle a geodesic makes with the meridians. 
We assume that $C$ contains a circle of latitude $\gamma_{0}$ where the function $r$ has a nondegenerate minimum. This means that $\gamma_0$ is a hyperbolic closed geodesic. We also assume that the other circles of latitude where $r$ has the same value $r_0$ as on $\gamma_0$ are not critical points of $r$. Thus $\gamma_0$ is the only closed geodesic on which $r = r_0$. 
The closed orbit of the geodesic flow of $C$ corresponding to the geodesic $\gamma_{0}$  has a homoclinic connection.
This means that the weak stable and weak unstable manifolds of this orbit coincide.
They lie in the set of vectors tangent to the geodesics for which the value of the Clairaut integral is $r_0$.
These geodesics are asymptotic to $\gamma_{0}$ as $t\rightarrow \pm\infty$.
One of them is shown in Figure 2.

Let $P$  and $Q$ denote the poles where the axis of revolution intersects $C$.
We parametrize $C$ with geodesic polar coordinates $(\theta,l)$, $\theta\in [0,2\pi)$ and $l\in [0,R]$. 
In these coordinates $C$ is determined by the profile function $r(l)$, which is required to be $C^{\infty}$, that gives the radius of a circle of latitude at a distance $l$ from the point $Q$ at the bottom of the surface.
For technical reasons we assume that there exists $l_1\in (0,R)$ and $b>0$, such that for $l\in (l_1-b,l_1+b)$ we have $r(l)\equiv r_1$, where $r_1$ is a constant such that $r_1>r_0$.
We will also assume that the value $l_0$ of $l$ on $\gamma_{0}$ satisfies $l_0<l_1-b$ and that $r$ is strictly increasing for $l\in (l_0,l_1-b)$.
Let $F$ denote the flat region given by those $p\in C$ such that $l(p)\in [l_1-b,l_1+b]$.
A typical shape for $C$ is shown in Figure 2.
\begin{figure}
\vspace{8cm}
\caption{Surface of revolution $C$.}
\end{figure}

Recall that the meridians of $C$ are all simple closed geodesics with the same period.

\begin{Lemma}There exist arbitrarily small smooth perturbations of the Riemannian metric of $C$ with support in $F$ such that:

(1) The meridians are preserved, i.e. all geodesics leaving from $P$ are simple, closed and with the same period, and coincide with the meridians outside $F$.

(2) $\gamma_{0}$ possesses a {\bf transverse} homoclinic orbit.\label{meridians}

\end{Lemma}

\PROOF The assumptions made above (in particular the requirement $r_1 > r_0$) and the properties of the Clairaut integral ensure that there are geodesics that pass through the flat region $F$ and are both forwards and backwards asymptotic to $\gamma_0$; indeed any vector in $F$ that makes angle~$\phi$ with the meridians will be tangent to such a geodesic if $r_1\sin\phi = \pm r_0$. Let $\gamma_{su}$ be such a geodesic and let $p$ be a simple point of $\gamma_{su}$ that lies in $F$. 
Parametrize $\gamma_{su}$ so that $\gamma_{su}(0) = p$.
Consider the frame $\{E_{1}, E_{2}\}$ in $T_{p}C$, where $E_{1}=\gamma_{su}'(0)$ and $E_{2}$ is the unit vector tangent to the meridian through $p$ pointing towards $P$.
Let $\{E_{1}(t), E_{2}(t)\}$ denote the frame in $T_{\gamma_{su}(t)}C$ (not orthonormal) obtained by parallel transport along $\gamma_{su}$.
Consider the map $f:{\bf R}^{2}\rightarrow C$ given by
\[f(t,x)=\exp_{\gamma_{su}(t)}(xE_{2}(t)).\]
For a small enough $\delta > 0$, the map $f$ is a diffeomorphism form $\Delta = (-\delta,\delta)^2$ to a neighborhood~$U \subset F$ of the point~$p$; the set~$U$ is shaded in Figure 2.
 For each fixed $t$, the curves $x\rightarrow f(t,x)$ are meridians parametrized by arc length.
In these coordinates the metric of $C$ satisfies:
\[g_{11}(t,x)=1,\;\;\;g_{12}(t,x)= a,\;\;\;g_{22}(t,x)=1,\]
where $a$ is a constant satisfying $-1 < a < 1$.

Let $\alpha : {\bf R}^2 \rightarrow {\bf R}$ be a smooth function with support inside $\Delta$
and let $g^\alpha$ be the metric defined by:
\[g^\alpha_{11}(t,x) = 1 - \alpha(t,x)x^{2},\;\;\;g^\alpha_{12}(t,x)= a,\;\;\;g^\alpha_{22}(t,x)=1.\]
 A simple computation shows that the Christoffel symbols of $g^\alpha$ satisfy:
 \[\Gamma^{1}_{22}\equiv \Gamma^{2}_{22}\equiv 0,\]
\[\Gamma^{1}_{11}(t,0)=\Gamma^{2}_{11}(t,0)=0.\]
Therefore the curve, $t\rightarrow f(t,0)$, and the curves, $x\rightarrow f(t_{0},x)$ for each fixed~$t_0$, are still geodesics, and thus the perturbation $g^\alpha$ preserves the geodesic $\gamma_{su}$ and the meridians.

A further computation shows that the Gaussian curvature of $g^\alpha$ at the point $\gamma_{su}(t) = f(t,0)$ is given by:
\begin{equation}
K_\alpha(t,0)=\frac{-\alpha(t,0)}{1-a^{2}}.  \label{gauss}
\end{equation}
This equation and Donnay's arguments from \cite{D} imply that for a suitable choice of the function~$\alpha$ the stable and unstable manifolds of $\gamma_{0}$ must have a {\sl transverse} intersection at the point $\gamma_{su}'(0)$, which concludes the proof of the lemma.

For the convenience of the reader, we give a brief sketch of Donnay's idea. 
Let $H^-$ and $H^+$ be the projections to $C$ of the strong stable and strong unstable manifolds associated to the geodesic $\gamma_{su}$.
Recall that the strong stable and strong unstable manifolds are given by the unit normal vectors to $H^-$ and $H^+$ that point to the same side as the tangent vector to $\gamma_{su}$.
 The geodesic curvatures at $\gamma_{su}(t)$ of $H^-$ and $H^+$ are solutions $u^-$ and $u^+$ of the Riccati equation
     \begin{equation}
u'(t) + u^2(t) + K(t) = 0, \label{riccati}
     \end{equation}
where $K(t)$ is the curvature at $\gamma_{su}(t)$. Before the perturbation we have $u^-_{old} \equiv u^+_{old}$. Let $t_1$ be the time when $\gamma_{su}$ enters and $t_2$ the time when $\gamma_{su}$ leaves the support of the perturbation of the metric. Then $u^-_{new}(t) =  u^-_{old}(t)$ for $t \geq t_2$ and  $u^+_{new}(t) =  u^+_{old}(t)$ for $t \leq t_1$. All that one needs to arrange is that $u^+_{new}(t_2) \neq u^-_{new}(t_2) = u^+_{old}(t_2)$. It is clear from (\ref{gauss}) and (\ref{riccati}) that this will be the case if the function $\alpha$ is chosen suitably.

\QED

\section{\em The main example}

As mentioned in the introduction, we begin with a one parameter family of surfaces of revolution $C_d$,  $0<d<\infty$, which have
the properties shown in Figure 1.
Each $C_d$ contains a region $R$ whose geometry is the same for all $d$. The boundary of $R$ is  a geodesic circle of latitude $\alpha$, which is shorter than any other geodesic circle of latitude in $R$: this ensures that any geodesic which enters $R$ must leave after a finite time. The region $R$ contains a flat cylinder $F$ and a hyperbolic closed geodesic $\gamma_{0}$; the distance of the surface from the axis of revolution is strictly decreasing as one moves from $F$ to $\gamma_{0}$. Attached to $R$ is a flat cylinder of length $d$ which is smoothly closed with a cap $D$; the geometry of $D$ is independent of $d$. The boundary of $D$ is a closed geodesic, which we shall denote by $\beta$. We choose $D$ so that $\beta$ is the only circle of latitude in $D$ that is a geodesic: this ensures that any geodesic which enters $D$ must leave after a finite time. Let $P$ be the center of the cap $D$ and $Q$ the center of the region $R$.

We parametrize $C_d$ with geodesic polar coordinates $(\theta,l)$ where $\theta \in {\bf R}/2\pi$ and $l$ is the geodesic distance of a point from $Q$. In these coordinates $C_d$ is determined by the profile function $r(l)$ which gives the distance from the axis of revolution of the circle of latitude which is at geodesic distance $l$ from $Q$. We  choose $r(l)$ to be $C^\infty$. Let  $l_\alpha$ be the value of $l$ on $\alpha$. Let $\rho = r(l_\alpha)$. Finally we introduce a third coordinate $\phi \in (-\pi/2,\pi/2]$ on the unit tangent bundle of $C_d \setminus \{P,Q\}$, which is the angle measured in the counterclockwise direction between a vector and the meridian  passing through the point where the vector is based (we shall take $\phi = \pi/2$ for the tangent vectors to the circle of latitude where this description is ambiguous).

Define $S_1$ and $S_4$ to be the sets of unit vectors with footpoints on $\alpha$ that point out of and into $R$ respectively. Define $S_3$ and $S_2$ to be the sets of unit vectors with footpoints on $\beta$ that point out of and into $D$ respectively. The coordinates $\theta$ and $\phi$ allow us to identify $S_i$ with $ ({\bf R}/2\pi) \times (-\pi/2,\pi/2)$ for $i = 1,\dots,4$.

For a vector $v \in S_1$, let us  follow the geodesic $\gamma_v$ defined by $v$ starting from time $t = 0$. We have chosen the geometry of $R$ and $D$ so that $\phi_t(v) = \gamma'_v(t)$ moves from $S_1$, through $S_2$, $S_3$ and $S_4$ in succession and then returns to $S_1$.
This process defines maps from 
$S_{i}$ to  $S_{i+1}$ for $i=1,\dots,4$~mod $4$,
which induce maps
\[\psi_{i}:({\bf R}/2\pi) \times (-\pi/2,\pi/2)\rightarrow ({\bf R}/2\pi) \times (-\pi/2,\pi/2),\;\;\;\;i=1,\dots,4.\]
Elementary trigonometry shows that
\begin{equation}
\psi_i(\theta,\phi) = (\theta + \frac{d\tan \phi}{\rho}, \phi) \qquad\hbox{for $i=1,3$.}     \label{cylinder}
\end{equation}

The properties of the Clairaut integral imply that
\begin{equation}
\psi_i(\theta,\phi) = (\theta + a_i(\phi), -\phi), \qquad\hbox{for $i=2,4$,} \label{minus}
\end{equation}    
where $a_2$ and $a_4$ are $C^\infty$ functions  satisfying 
\begin{equation}
a_2(0)= \pi = a_4(0).
\label{aa}
\end{equation}
 Notice that $a_2$ and $a_4$ {\bf do not} depend on $d$.

As explained in the introduction, we make a small perturbation of the surfaces of revolution $C_d$.
Recall from Section 3 (cf. Lemma \ref{meridians}) that there exist arbitrarily small, smooth perturbations  of the metric, such that:

(1) The support of the perturbation is confined to the interior of the flat cylinder $F \subset R$.

(2) The meridians are preserved, i.e. all geodesics leaving from $P$ are simple, closed and with the same period, and coincide with the meridians outside $F$.

(3) $\gamma_{0}$ possesses a transverse homoclinic orbit.

We can also assume that the perturbed metric has the following property:

 (4) There is a constant $\phi_0 > 0$ such that every geodesic that enters $R$ through $\alpha$ with $|\phi| < \phi_0$ must exit $R$ after a finite time.

To see this, note that we can choose $\phi_0 > 0$, $T_0 > 0$ and $\ep_0 > 0$ so that every geodesic of the original metric that enters $R$ with $|\phi| \leq 2\phi_0$ must exit the $2\ep_0$ neighborhood of $R$ by time $T_0$ at the latest. If the perturbation is small enough, every geodesic of the new metric that crosses $\alpha$ with $|\phi| \leq \phi_0$ must exit the  $\ep_0$ neighborhood of $R$ by time $2T_0$ at the latest.

The above properties ensure that the orbits of the geodesic flow that start in $S_1$ with $|\phi| < \phi_0$ still pass through $S_2$, $S_3$ and $S_4$ in succession and then return to $S_1$. The transitions from $S_i$ to $S_{i+1}$, for $i = 1,2,3$, are exactly the same as for the unperturbed metric and do not change $|\phi|$. Property (4) means that every orbit that leaves $S_4$ with $|\phi| < \phi_0$ must pass through $S_1$.  Let $\tilde{\psi}_{4}:({\bf R}/2\pi) \times (-\phi_0,\phi_0)\rightarrow ({\bf R}/2\pi) \times (-\pi/2,\pi/2)$ be the map induced by the transition from $S_{4}$ to $S_{1}$ for the perturbed surface.
Notice that $\tilde{\psi}_{4}$ is determined by the metric in $R$ and is independent of $d$.

The next result will be crucial for our construction:

\begin{Lemma}There exists $d_{0}>0$ and a perturbation of the metric in $R$ satisfying properties (1)--(4) above, such that for any $d > d_{0}$, the map
\[ \tilde\Psi_d \df \tilde \psi_4 \circ \psi_3 \circ \psi_2 \circ \psi_1 : ({\bf R}/2\pi) \times (-\phi_0,\phi_0) \to ({\bf R}/2\pi) \times (-\pi/2,\pi/2)\]
has a homotopically nontrivial invariant circle in the set 
\( ({\bf R}/2\pi) \times  [-1/d,1/d] . \) \label{key}
\end{Lemma}

\PROOF
For $d>0$ we define a scaling map $\varphi_{d} : ({\bf R}/2\pi) \times {\bf R} \to ({\bf R}/2\pi) \times {\bf R}$ by $\varphi_{d}(\theta,\phi) \df (\theta,d\phi)$.
Observe that $\varphi_{d}$ maps the region $| \phi| \leq 1/d$ diffeomorphically onto the set \[ S\df ({\bf R}/2\pi) \times [-1,1].\] In order to show that $\tilde\Psi_d$ has an invariant circle in the region $|\phi| < 1/d$, it suffices to show that the map \[\tilde k_d \df \varphi_d \circ \tilde\Psi_d \circ \varphi_d^{-1} : S \to ({\bf R}/2\pi) \times {\bf R},\] which is well defined for $d > 1/\phi_0$, has an invariant circle. 

Let $\Psi_d =  \psi_4 \circ \psi_3 \circ \psi_2 \circ \psi_1$ and let $k_d = \varphi_d\circ \Psi_d \circ \varphi_d^{-1}$. 

Our first step is to study the limiting behaviour of $k_d$ and $\tilde k_d$ as $d \to \infty$ using:

\begin{Lemma}Let $f=(g,h):({\bf R}/2\pi) \times (-\phi_0,\phi_0)\rightarrow ({\bf R}/2\pi) \times (-\pi/2,\pi/2)$ be a $C^{\infty}$ map with the property
\[h(\theta,0)=0,\]
for all $\theta\in {\bf R}/2\pi $.
For $d>1/\phi_0$ define
\[f_{d}=\varphi_{d}\circ f\circ\varphi_{d}^{-1}:S\rightarrow ({\bf R}/2\pi) \times {\bf R}.\]
Then $f_{d}$ converges in the $C^{\infty}$ topology as $d \to \infty$ to the map
\[(\theta,\phi)\rightarrow (g(\theta,0),\phi\frac{\partial h}{\partial \phi}(\theta,0)).\]    \label{technical}
\end{Lemma}

\PROOF
Since $h(\theta,0) = 0$, there is a $C^\infty$ function $H(\theta,\phi)$ such that
\begin{equation}
 h(\theta,\phi) = \phi H(\theta,\phi).
\label{key*}
\end{equation}
We have
\begin{equation}
f_d(\theta,\phi) = (g(\theta,\phi/d),dh(\theta,\phi/d)).
\label{key**}
\end{equation}
It is easy to see that as $d \to \infty$ the first component of $f_d$ converges in the $C^\infty$-topology to the function
\[  (\theta,\phi) \to g(\theta,0).\]
It follows from (\ref{key*}) and (\ref{key**}) that the second component of $f_d$ is
\[   (\theta,\phi) \to \phi H(\theta,\phi/d),\]
which converges in the $C^\infty$-topology to the function
\[   (\theta,\phi) \to \phi H(\theta,0) = \phi \frac{\partial h}{\partial \phi}(\theta,0).\]
\QED

Now observe that $k_{d}$ is the composition of the four maps  $\varphi_{d}\circ\psi_{i}\circ\varphi^{-1}_{d}$, $i=1,\dots,4$.
It is easily seen from (\ref{cylinder}) that for $i = 1$ and $3$
\[\varphi_{d}\circ\psi_{i}\circ\varphi^{-1}_{d}(\theta,\phi)=(\theta+\frac{d}{\rho}\tan\frac{\phi}{d},\phi),\]
which converges in the $C^{\infty}$ topology to the map
\[(\theta,\phi)\rightarrow (\theta+\frac{\phi}{\rho},\phi).\]

It follows from equations (\ref{minus}), (\ref{aa}) and Lemma \ref{technical} that for $i = 2$ and $4$ the map $\varphi_{d}\circ\psi_{i}\circ\varphi^{-1}_{d}$ converges to the map
\begin{equation}
(\theta,\phi)\to (\theta+\pi,-\phi).
\label{eins}
\end{equation}

A simple calculation now shows that as $d \to \infty$ the map $k_d$ converges in the $C^\infty$ topology to the map
\[ k_\infty  (\theta,\phi) = ( \theta + \frac{2\phi}{\rho} , \phi ). \]

Observe that $k_\infty$ is a nondegenerate twist map. By Moser's twist theorem \cite{Mo}, there is a neighbourhood $\cal U$ of $k_\infty$ in the $C^\infty$ topology such that every map in $\cal U$ possesses a homotopically nontrivial invariant circle.

Now consider the perturbed surface. 
Let $\tilde{\psi}_{4}=(\tilde{g},\tilde{h})$.
It follows immediately from Lemma \ref{technical} that $\varphi_{d}\circ\tilde{\psi_{4}}\circ\varphi^{-1}_{d}$ converges to the map
\begin{equation}
(\theta,\phi)\to (\tilde{g}(\theta,0),\phi\frac{\partial \tilde{h}}{\partial \phi}(\theta,0)).
\label{zwei}
\end{equation}

We now know that $\tilde k_d$ is the composition of $\varphi_{d}\circ\tilde\psi_{4}\circ\varphi^{-1}_{d}$ and $\varphi_{d}\circ\psi_{i}\circ\varphi^{-1}_{d}$, $i=1,2,3$, and each of these maps approaches a limit as $d \to \infty$. Hence there is a map $\tilde k_\infty$ such that $\tilde k_d$ converges to $\tilde k_\infty$ in the $C^\infty$ topology as $d \to \infty$. If the perturbation of the metric is  small enough, $\tilde\psi_4$ will be close  to $\psi_4$,  the map in (\ref{zwei}) will be close to the map in (\ref{eins}), and $\tilde k_\infty$ will be in $\cal U$. It follows that for such a perturbation, $\tilde k_d$ belongs to $\cal U$, and thus has a homotopically nontrivial invariant circle, for all large enough $d$. This completes the proof of Lemma \ref{key}.

\QED

{From} now on we shall assume that $d \geq d_0$ so the conclusion of Lemma \ref{key} holds.
Let $\tilde{C}_{d}$ denote the perturbation of $C_{d}$ we constructed above.
Observe that the invariant circle of the map $\tilde\Psi_d$ in the region $|\phi | \leq 1/d$, gives rise to an invariant torus $T_{d}$ for the geodesic flow of $\tilde{C}_{d}$. 
The projection map $\pi\mid_{T_{d}} : T_{d} \to \tilde{C}_{d}$ becomes singular on two simple closed curves which project to simple closed curves $\Gamma_{P,d}$ and $\Gamma_{Q,d}$ as indicated in Figure 1.
If $d$ is large enough, $\Gamma_{P,d}$ will lie in $D$ and $\Gamma_{Q,d}$ will lie in $R$.
Henceforth we assume that all $d\geq d_{0}$ have this property.
We shall also drop the subscript $``d"$ in order to simplify the notation.

Let $\Theta: S\tilde{C}_{d} \to S\tilde{C}_{d}$ denote the flip:
\[\Theta(p,v)=(p,-v).\]
Then $\Theta T_{d}$ is also an invariant torus of the geodesic flow, and the projection map $\pi\mid_{\Theta T_{d}}$ becomes singular above the curves $\Gamma_{P}$ and $ \Gamma_{Q}$.

Let $V_{P}$ and $V_{Q}$ be the open neighborhoods of $P$ and $Q$ that are bounded by $\Gamma_{P}$ and $\Gamma_{Q}$ respectively.
Clearly, $T_{d}$ and $\Theta T_{d}$ separate the unit sphere bundle of 
$\tilde{C}_{d}$ into two invariant sets, $W_{1}$ and $W_{2}$. 
One of them, which we shall call $W_{1}$, contains $S(V_{P}\cup V_{Q})$.
If $v$ is in $W_{1}$ the geodesic $\gamma_{v}$ oscillates between $D$ and $R$; in particular $\gamma_{v}$ enters both $D$ and $R$.
It is clear from this that the set $W_{2}$ contains the unit vectors tangent to $\gamma_{0}$, and thus the transverse homoclinic orbit associated with $\gamma_{0}$.

\begin{Lemma} There exists $t_{0}>0$, such that if $\gamma$ is any geodesic leaving from a point in $V_{P}\cup V_{Q}$, then $\gamma$ spends time at most $t_{0}$ inside $D\cup R$ (in each passage). The time $t_0$ is independent of $d$. \label{time}
\end{Lemma}
  
\PROOF The tangent vectors to $\gamma$ belong to the invariant set $W_1$ defined above.
Observe that $\overline{W}_1$ is compact and does not contain the tangent vectors to $\alpha$ and $\beta$, the boundaries of $R$ and $D$ respectively.

Both regions $R$ and $D$ have the property that their boundary is a closed geodesic and that any geodesic which crosses the boundary and enters the region must exit after finite time.
Since geodesics always cross transversally, it is easy to see that the time between entry and exit is a continuous function of the tangent vector at the time of entry.
Let $t_{0}$ be the supremum of the longest geodesic segment in $D\cup R$ that is part of a geodesic whose tangent vectors are in $\overline{W}_1$.
It is easy to see from the above that $t_{0}<\infty$ and each passage of $\gamma$ through $D\cup R$ takes time at most $t_{0}$.

\QED

\begin{Lemma}Given $\varepsilon>0$, there exists $d(\varepsilon)$ such that for $d>d(\varepsilon)$, the geodesic flow of $\tilde{C}_{d}$ satisfies:
\[h_{top}({W_{1}})<\varepsilon.\]  \label{entropy}
\end{Lemma}

\PROOF On account of the variational principle for entropy and the relationship between Liapunov exponents and entropy, it suffices to show that given $\varepsilon>0$, there exists $d(\varepsilon)$ such that for $d>d(\varepsilon)$ and all $v\in W_{1}$ we have:
\[\limsup_{t\rightarrow +\infty}\frac{1}{t}\log \parallel d\phi_{t}(v)\parallel<\varepsilon.\]

We follow an argument of Manning \cite{Manning}.
Let $K$ denote the Gaussian curvature of $\tilde{C}_{d}$ and let $Y$ denote any Jacobi field along the geodesic $\gamma_{v}$ defined by $v$.
Also let $L$ be an upper bound for $ K^2$ in the region $R\cup D$. We may assume that $L > 1$ and $\varepsilon < 1$.

Consider the function:
\[ y_{\varepsilon}(t) \df \varepsilon^{2} \langle  Y(t),Y(t) \rangle  + \langle  Y'(t),Y'(t) \rangle. \]
Using the Jacobi equation we deduce:
\[ y_{\varepsilon}'(t) = 2\varepsilon^{2} \langle  Y'(t),Y(t)\rangle  + 2\langle  -K(t)Y(t),Y'(t) \rangle. \]
Let  $t$ be such  that $\gamma_{v}(t)\in R\cup D$. Since $\varepsilon < L$, we have
\[ |y_{\varepsilon}'(t)| \leq 2\varepsilon^{2} \| Y'(t)\|\| Y(t)\| + 2L^{2}\| Y'(t)\|\| Y(t)\| \leq 4L^{2}\| Y'(t)\|\| Y(t)\| \]
\[ = \frac{2L^{2}}{\varepsilon} (2\varepsilon \| Y'(t)\|\| Y(t)\|) \leq \frac{2L^{2}}{\varepsilon} (\varepsilon^{2} \langle  Y(t),Y(t)\rangle  + \langle  Y'(t),Y'(t) \rangle) \]
\begin{equation}
= \frac{2L^{2}}{\varepsilon}y_{\varepsilon}(t).  \label{cinco}
\end{equation}

On the other hand, if $t$ is such that $\gamma_{v}(t)$ belongs to the flat cylinder, we have:
\[  y_{\varepsilon}'(t) = 2\varepsilon^{2}\langle  Y'(t),Y(t) \rangle  , \]
and thus
\begin{equation}
 |y_{\varepsilon}'(t)| \leq 2\varepsilon^{2}\| Y'(t)\|\| Y(t)\| 
\leq \varepsilon(\varepsilon^{2} \langle  Y(t),Y(t) \rangle  + \langle  Y'(t),Y'(t) \rangle) = \varepsilon y_{\varepsilon}(t).  \label{seis}
\end{equation}
Next observe that
\[ \limsup_{t\rightarrow + \infty}\frac{1}{t}\log \| d\phi_{t}(v)\| \leq \limsup_{t \to +\infty} \frac{1}{t} \log\sqrt{\frac{y_{\varepsilon}(t)}{\varepsilon^{2}}}\]
\begin{equation}
=\limsup_{t\rightarrow+\infty}\frac{1}{2t}\log y_{\varepsilon}(t)=\limsup_{t\rightarrow+\infty}\frac{1}{2t}\int_{0}^{t}\frac{y_{\varepsilon}'(s)}{y_{\varepsilon}(s)}\,ds.  \label{siete}
\end{equation}

Recall that the geodesic $\gamma_{v}$ oscillates between $D$ and $R$.
In each passage through the flat cylinder between $D$ and $R$, $\gamma_{v}$ spends time at least $d$ and by Lemma \ref{time}, it spends time at most $t_{0}$ during each visit to $R\cup D$. 
By combining equations (\ref{cinco}), (\ref{seis}) and (\ref{siete}) we deduce:
\[\limsup_{t \to +\infty}\frac{1}{t}\log \| d\phi_{t}(v) \| \leq \frac{L^{2}t_{0}}{d\varepsilon} + \frac{\varepsilon}{2}.\]
If we take $d > d(\varepsilon) \df 2L^{2}t_{0} / \varepsilon^{2}$, we obtain
\[ \limsup_{t \to +\infty}\frac{1}{t}\log \| d\phi_{t}(v)\| < \varepsilon. \]

\QED

We are ready to prove the main result of this section:

\begin{Theorem}For all $d$ sufficiently large and $p\in V_{P}\cup V_{Q}$,
\[\limsup_{T\rightarrow\infty}\frac{1}{T}\log n_{T}(p,q)\leq 
\limsup_{T\rightarrow\infty}\frac{1}{T}\log\int_{\tilde{C}_{d}}n_{T}(p,q)\,dq< h_{top},\]
for a.e. $q\in \tilde{C}_{d}$.
\end{Theorem}

\PROOF Since $\pi^{-1}(V_{P}\cup V_{Q}) \subset W_1$, we see from  Proposition \ref{p1} that it suffices to show that for all $d$ sufficiently large,
\[h_{top}({W_{1}})<h_{top}.\]

Let $h$ denote the entropy of the horseshoe associated with $\gamma_{0}$; obviously $h$ is independent of $d$ and $h_{top}(\phi_{t})\geq h$. 
The theorem now follows from the last lemma.

\QED

This theorem gives negative answers to Questions I, I$'$ and I$''$ in the Introduction.

Now observe that the tori $T_d$ and $\Theta T_d$ that were considered above will survive under any sufficiently small (in the $C^\infty$ topology) perturbation of the metric that we constructed above. These tori will still separate the unit tangent bundle and the conclusions of Lemma 4.4 and Theorem 4.5 carry over. It is clear from this that Question II in the introduction has a negative answer.

Finally we observe that since $h_{top}(S_p) \leq h_{top}(W_1)$ for any $p \in V_p \cup V_q$, we have $h_{top}(S_p) < h_{top}$ for any $p \in V_p \cup V_q$.

\vspace{1cm}
\noindent Keith Burns\\
Mathematics Department\\
Northwestern University\\
Evanston IL 60208\\
U.S.A.\\
E-mail: burns@math.nwu.edu

\vspace{1cm}

\noindent Gabriel P. Paternain\\
IMERL-Facultad de Ingenier\'\i a\\
Julio Herrera y Reissig 565, C.C. 30\\
Montevideo\\
Uruguay\\
E-mail: gabriel@cmat.edu.uy

\end{document}